\documentclass{amsart}

\usepackage{amssymb}

\usepackage{graphicx}

\newtheorem{theorem}{Theorem}[section]
\newtheorem{lemma}[theorem]{Lemma}

\theoremstyle{definition}

\newtheorem{example}[theorem]{Example}

\theoremstyle{remark}
\newtheorem{remark}[theorem]{Remark}

\theoremstyle{problem}

\def\CC{\mathbb{C}}

\def\PP{\mathbb{P}}

\def\CP{\mathbb{CP}}

\def\dd{\partial}

\def\cal{\mathcal}  
\newcommand{\A}{{\cal A}}

\numberwithin{equation}{section}

\newcommand{\codim}{\operatorname{codim}}

\newcommand{\Proj}{\operatorname{Proj}}
\newcommand{\Sat}{\operatorname{Sat}}

\newcommand{\wt}{\widetilde}
\newcommand{\ol}{\overline}
\newcommand{\PGL}{\operatorname{PGL}}

\newcommand{\HP}{\operatorname{HP}}

\begin{document}

\title[The Jacobian Ideal]{The Jacobian ideal of a hyperplane arrangement}


\author{Max Wakefield}
\address{Department of Mathematics, Hokkaido University, Sapporo 060-0810, Japan}
\curraddr{}
\email{wakefield@math.sci.hokudai.ac.jp}

\author{Masahiko Yoshinaga}
\address{The Abdus Salam ICTP, Strada Costiera 11, Trieste 34014, Italy}
\curraddr{}
\email{myoshina@ictp.it}

\thanks{The first author has been supported by NSF grant \# 0600893
  and the NSF Japan program. The second author has been supported 
by JSPS Postdoctoral Fellowships for Research Abroad. }

\subjclass[]{}

\date{}

\dedicatory{}

\begin{abstract}
The Jacobian ideal of a hyperplane arrangement is an ideal 
in the polynomial ring whose generators are the partial 
derivatives of the arrangements defining polynomial. 
In this article, we prove that an arrangement can be 
reconstructed from its Jacobian ideal. 
\end{abstract}

\maketitle


\section{Introduction}

Let $V\cong \CC^\ell$ and choose coordinates for $V^*$ such that we
can identify the symmetric algebra $S=S(V^* )$ with the polynomial
ring $\CC [z_1,\ldots ,z_\ell]$. A hyperplane in $V$ is a codimension
one affine space in $V$. A hyperplane arrangement in $V$ is a finite
collection of hyperplanes denoted by $\A$. When all the
hyperplanes of an arrangement contain the origin we say the
arrangement is central. For most of this note we assume the
arrangement is central. In this case we can `projectivize' all the
hyperplanes and view the arrangement as an arrangement of hyperplanes
in $\CP^{\ell -1}$. Further, we say a central arrangement $\A$ is essential if $\bigcap\limits_{H\in \A}H=\{ 0\}$. 

For each $H\in \A$ choose a linear polynomial
$\alpha_H\in S$ such that $H=\ker \alpha_H$. Let $Q=\prod\limits_{H\in
  \A}\alpha_H$ denote the defining polynomial of the arrangement
$\A$. Then the main character of this note is the homogeneous ideal in
$S$ defined by $$J(Q):=\left( \frac{\dd Q}{\dd z_1},\ldots,\frac{\dd
    Q}{\dd z_\ell}\right)_.$$ We call this ideal the Jacobian ideal and sometimes denote it by $J(\A )$. 
The Jacobian ideal determines a closed subscheme 
$\Proj S/J(Q)$ of the projective space $\CP^{\ell-1}$, which we 
call the Jacobian scheme.

The purpose of this paper is to prove the following 
result, which simply put, states that the Jacobian scheme 
contains all the information of the arrangement. We say two hyperplane arrangements $\A_1$ and $\A_2$ are \emph{identical} when $Q_1=cQ_2$ for some $c\in \CC^*$ where $Q_1$ and $Q_2$ are the defining polynomials respectively. 

\begin{theorem}\label{main}

Suppose $\A_1$ and $\A_2$ are two central and essential arrangements in dimension
$\ell \geq 3$. Then $\A_1$ and $\A_2$ are identical if and only if the
Jacobian schemes $\Proj S/J(\A_1)$ and $\Proj S/J(\A_2)$ are equal 
as closed subschemes of $\CP^{\ell-1}$. 

\end{theorem}

The proof of Theorem \ref{main} is inspired by a Torelli-type 
theorem of Dolgachev and Kapranov \cite{DK, D}. 
In \cite{DK}, Dolgachev and Kapranov prove that the module 
$D(\A)$ of 
derivations of a generic arrangement $\A$ contains all 
the information of 
the arrangement. More precisely, they consider the 
set of jumping lines of 
the torsion free (actually locally free when $\A$ is 
a generic arrangement \cite{mus-sch, Yuz}) sheaf 
$\wt{D}(\A)$ on the projective space. 
From the set of jumping lines, the arrangement $\A$ can be 
recovered. Then in \cite{D}, by considering a certain subsheaf of
$\wt{D}(\A)$,  Dolgachev extended these results to a wider class of arrangements.

Instead of jumping lines, we consider the subscheme obtained as the intersection 
$\overline{K}\cap\Proj S/J(\A)\subset\CP^{\ell-1}$,  for a given hyperplane 
$K\subset V$. 
In particular, we focus on the $(\ell-3)$-dimensional 
components of $\overline{K}\cap\Proj S/J(\A)$ and compute 
the $(\ell-3)$-dimensional degree (i.e. the coefficient of the $\ell -3$ term in the Hilbert polynomial). Then we can prove that $K\in\A$ precisely 
when this degree is maximized. We also note that 
the reduced Jacobian scheme $\Proj S/\sqrt{J(Q)}$ 
does not contain all the information of $\A$ (see Remark \ref{rem:rad}). 

Another closely related result is found in \cite[Prop. 1.1]{Dn}. 
Let $f\in S_d$ be a homogeneous polynomial of degree $d$. 
Then Donagi proved that the Jacobian ideal $J(f)$ recovers $f$ 
up to $\PGL$-action. Our main result in this paper 
strengthens this assertion for hyperplane arrangements, 
namely, the saturated Jacobian ideal $\Sat(J(Q))$ 
recovers the defining equation $Q$ up to constant 
multiple.

\medskip

At this time the authors would like to thank T. Abe, H. Schenck, B. Shelton, H. Terao, K. Ueda and S. Yuzvinsky for many helpful discussions. The authors would also like to thank the creators of the computer algebra system Macaulay 2 (see \cite{GS}) since many of the ideas from this note originated by computing examples in this program.

\section{Minimal Components of $J(Q)$}\label{PrimaryJ(Q)}

In this section, we will study the minimal primary 
components of the Jacobian ideal of the arrangement $\A$. 

Throughout this paper we use the following notation. Let $L(\A)$ be the intersection lattice of $\A$ which is the set of all intersections of elements from $\A$ with the order being reverse inclusion. Moreover, let $L_k(\A)=\{ X\in L(\A)\mid  \codim (X) =k\}$. For $X\in L(\A)$ let $\A_X=\{H\in \A\mid  X\subseteq H\}$ and $L(\A)_X=\{ Y\in L(\A)\mid X\subsetneq Y\}$. Then we define the M\"obius function $\mu$ on $L(\A)$ by setting $\mu (V)=1$ and the recursive formula: $$\mu (X)=-\sum\limits_{Y\in L(\A)_X}\mu (Y).$$


We assume the dimension $\ell\geq 3$. For given an intersection $X\in L(\A)$, put 
\begin{eqnarray*}
Q_X&=&\prod_{H\in\A_X} \alpha_H, \\
\ol{Q}_X&=&\frac{Q}{Q_X}=\prod_{X\nsubseteq H}\alpha_H. 
\end{eqnarray*}
Obviously $Q=Q_X \ol{Q}_X$. 
Let us denote $I(X):=\sum_{H\in\A_X}S\alpha_H$ 
the prime ideal representing $X$. 

Since the Jacobian ideal $J(Q)$ determines the 
singular loci of the union $\bigcup_{H\in\A} H$ of 
hyperplanes, we have 
$$
\sqrt{J(Q)}=\bigcap_{X\in L_2(\A)}I(X). 
$$
This implies that the set of minimal associated 
primes of $J(Q)$ is $\{ I(X)\mid X\in L_2(\A)\}$. 
The localization technique enables us to obtain the corresponding minimal primary components as follows. 
\begin{lemma}
\label{mincomp}
The set of minimal components of the Jacobian ideal $J(Q)$ is equal to $\{J(Q_X)\mid  X\in L_2(\A)\}$. 
\end{lemma}

\begin{remark}
Generally, the Jacobian ideal $J(Q)$ has a lot of 
embedded primes. If $\A$ is a free arrangement, then 
$S/J(Q)$ is known to be Cohen-Macaulay \cite{T}. In this 
case, $J(Q)$ has no embedded primes. Thus we have 
the primary decomposition $J(Q)=\bigcap_{X\in L_2(\A)}J(Q(X))$, 
\cite{T}. 
\end{remark}

\begin{remark}
The degree of the ideal $J(Q_X)$ is 
$$
\deg J(Q_X)=\mu(X)^2=(|\A_X|-1)^2. 
$$
Hence the degree of the Jacobian ideal $J(Q)$ is 
$\sum_{X\in L_2(\A)}\mu(X)^2$. For details see 
\cite[Theorem 2.5]{S-elem}. 
\end{remark}


\section{$\Proj S/J(Q)$ intersected with a hyperplane}

Fix a hyperplane $K=\{ \beta=0\}$ that is not necessarily in $\A$. 
In this section, we consider the codimension two components 
of $\ol{K}\cap\Proj S/J(Q)=\Proj S/(J(Q)+(\beta))$ 
in $\CP^{\ell-1}$. In particular, we compute its degree 
in terms of the M\"obius function. 

The essential part of the computation is the 
following $2$-dimensional case. 

\begin{lemma}
\label{lem:2dim}
Let $Q(z_1, z_2)=a_0z_1^n+a_1z_1^{n-1}z_2+\ldots +a_nz_2^n
\in\CC[z_1, z_2]$ be a non-zero degree $n$ 
homogeneous polynomial of two variables. Suppose 
$\{Q=0\}$ defines a distinct $n$ lines. 
Then $J(Q)+(z_2)=(z_1^{n-1}, z_2)$ and 
$$
\dim \CC[z_1, z_2]/(J(Q)+(z_2))=n-1. 
$$
\end{lemma}


Recall that if $I\subset S$ is a homogeneous ideal and 
assume $\dim \Proj S/I\leq m$, then the 
Hilbert polynomial is of the form 
$$
\HP(S/I, d)=\frac{a_m}{m!}d^m+
\frac{a_{m-1}}{(m-1)!}d^{m-1}+\cdots. 
$$
Let us denote the coefficient $a_m$ by 
$\deg_{(m)}\Proj S/I$, which depends only on the 
$m$-dimensional components of the closed subscheme 
$\Proj S/I\subset \CP^{\ell-1}$. By definition, 
if $\dim \Proj S/I< m$, then $\deg_{(m)}\Proj S/I=0$.

\begin{lemma}\label{degree}

For any arrangement $\A$ with defining polynomial $Q$ and any  hyperplane $K=\{\beta=0\}$, we have 
\begin{equation}\label{eq:mu}
\deg_{(\ell-3)} \Proj S/(J(Q)+(\beta)) =
\sum_{X\in L_2, X\subset K}\mu(X).
\end{equation}

\end{lemma}

\proof First note that every $(\ell-3)$-dimensional component 
of $\Proj S/(J(Q)+(\beta))$ is 
of the form $\Proj S/(J(Q_X)+(\beta))\subset\CP^{\ell-1}$ 
such that $X\in L_2(\A)$ and $X\subset K$. Then the lemma is 
immediate from Lemma \ref{lem:2dim}. \qed

We denote the right hand side of (\ref{eq:mu}) by 
$\mu_\A(K):=\sum_{X\in L_2, X\subset K}\mu(X)$.

\section{Reconstruction of the arrangement by the Jacobian scheme}

In this section we prove that the Jacobian ideal of a hyperplane arrangement and its saturation contain all the information from the arrangement. Hence, we prove Theorem \ref{main}. Let $\A$ be a central arrangement.  

\begin{lemma}\label{muLL}

If the hyperplane $K$ is in $\A$, then $\mu_{\A}(K)=|\A |-1$. If $K$ is not in $\A$, $\ell \geq 3$ and $\A$ is essential, then $\mu_{\A}(K)<|\A |-1$.

\end{lemma}

\proof 
The first statement 
follows easily from the definition (\ref{eq:mu}). 
Suppose that $K$ is not in $\A$. Put the 
set $L_2(\A)^K=\{X\in L_2(\A)\mid X\subset K\}$. 
If $L_2(\A)^K$ is empty, there is nothing to prove. 
If $L_2(\A)^K=\{X\}$ consists 
of one element, then there exists $H\in\A$ such that 
$X\nsubseteq H$ since $\A$ is essential. Hence $|\A_X|\leq |\A|-1$. 
We also obtain $\mu_\A(K)=|\A_X|-1<|\A|-1$. 
Finally suppose $L_2(\A)^K=\{X_1, X_2, \ldots, X_p\}$ with 
$p\geq 2$. 
Then from the assumption, we have $\A_{X_i}\cap\A_{X_j}=\emptyset$ 
for $1\leq i<j\leq p$. Thus we have 
$\mu_\A(K)=\sum_{i=1}^p|\A_{X_i}|-p<|\A|-1$. \qed

Now, we can prove Theorem \ref{main}. Let $\A$ be an 
essential hyperplane arrangement with $\ell\geq 3$. 
Let $K=\{\beta=0\}$ be a hyperplane and $\ol{K}\subset\CP^{\ell-1}$ 
the projectivization. Then the scheme theoretic 
intersection with $\Proj S/J(Q)$ is obtained by 
$$
\ol{K}\cap\Proj S/J(Q)=\Proj S/(J(Q)+(\beta)). 
$$
From Lemma \ref{muLL}, $\deg_{(\ell-3)}\ol{K}\cap\Proj S/J(Q)$ 
is not greater than $|\A|-1$ and maximized precisely when 
$K\in\A$. This reconstructs $\A$ from $\Proj S/J(Q)$. 
\qed

\begin{example}\label{rem:rad} It may be worth noting that from the reduced 
Jacobian scheme $\Proj S/\sqrt{J(Q)}$, 
we can not reconstruct $\A$. 
Suppose $\A_1$ is defined by $Q_1=z_1 z_2 z_3 (z_1+z_2-z_3)$ and 
$\A_2$ is defined by $Q_2=Q_1\times (z_1-z_3)$. 
Recall in general $\Proj S/\sqrt{J(Q)}$ is the reduced scheme 
structure on the singular locus, which is the union 
of codimension two intersections $X\in L_2(\A)$. 
Then the radical of Jacobian ideals are equal, 
more precisely, 
$$
\sqrt{J(Q_1)}=
\sqrt{J(Q_2)}=
(z_1,z_2)\cap (z_1,z_2-z_3)\cap (z_2,z_1-z_3)\cap (z_1,z_3)\cap (z_2,z_3)\cap (z_1+z_2,z_3). 
$$ So, the reduced Jacobian ideal does not even record the number of hyperplanes.
\end{example}

\begin{example}
Let $\A_1$ and $\A_2$ be arrangements of generic five planes 
in $\CC^3$. We may assume that $\A_1$ and $\A_2$ are not 
projectively equivalent (since $\dim\PGL(3,\CC)=8$ is less than 
$10=$ 
the dimension of the configuration space of five planes). 
On the other hand, the scheme $\Proj S/J(\A_i)$ is just 
ten points with the constant structure sheaf. Hence 
$\Proj S/J(\A_1)\simeq \Proj S/J(\A_2)$ as schemes. 
The authors do not know whether if there exist such 
pairs in higher dimension. 
\end{example}

\bibliographystyle{amsplain}

\end{document}